\documentclass[12pt]{amsart}
\usepackage{amsmath,amssymb,amsthm,url}
\def\C{{\mathbb C}}

\long\def\comment#1\endcomment{}

    \theoremstyle{theorem}
         \newtheorem{theorem}{Theorem}

    \theoremstyle{definition}
         \newtheorem{remark}[theorem]{Remark}

\begin{document}

\title{On some results of S. Abramyan and T. Panov}

\author{A. Skopenkov}

\thanks{Homepage: \texttt{https://users.mccme.ru/skopenko}.
Moscow Institute of Physics and Technology, Independent University of Moscow.
Supported in part by the Russian Foundation for Basic Research Grant No. 19-01-00169 and by Simons-IUM Fellowship.}

\date{}

\maketitle

\abstract
This note is purely expository and is an extended version of  math review to the paper [AP19]=arXiv:1901.07918v3 by S. Abramyan and T. Panov published in Proc. of Steklov Math. Inst. 305 (2019).
The authors construct simplicial complexes for whose moment-angle complexes certain homotopy classes are non-trivial.
I present in a shorter and clearer way the main definition and the statement of Theorem 5.1 from [AP19].
The clarification reveals that the main definition used in the statements
of the main results is not given [AP19].
\endabstract

\begin{remark}[Motivations]\label{r:mot}
This note is purely expository and is an extended version of math review to \cite{AP19}.
I bear no responsibility for results of \cite{AP19}.
I present in a shorter and clearer way the main definition and the statement of Theorem 5.1 from \cite{AP19}.
The clarification reveals that the main definition used in the statements
of the main results is not given \cite{AP19} (see justification in Remark \ref{r:spec}).  
Since the authors refused to update \cite{AP19},\footnote{Although it is not a task of math reviewer, I spent quite some efforts on persuading the authors to publish such an update.}
this note serves as an invitation to publish reliable and clear statements and proofs of the results whose non-rigorous statements are given in \cite{AP19}
(of course properly mentioning contribution by S. Abramyan and T. Panov).
One of the best estimations of how significant flaws are is the amount of time required for authors
(or for others) to make a corrected version publicly available upon request of a reviewer.
\end{remark}

The main results of \cite{AP19} (Theorem 5.1, 5.2 and 7.1) are constructions of simplicial complexes $K$ for whose moment-angle complexes $\mathcal Z_K$ (defined just below) certain homotopy classes are non-trivial.

Let $K$ be a simplicial complex, $V(K)$ its vertex set, $v$ its vertex and $\sigma$ its face.
Denote
$$
\Delta_{\sigma,v}:=\begin{cases} D^2 & v\in\sigma \\  S^1 & v\not\in\sigma
\end{cases}
\quad\text{and}\quad
\mathcal Z_K:=\bigcup\limits_{\sigma\in K}\prod_{v\in V(K)}\Delta_{\sigma,v}\subset (D^2)^{|V(K)|}.
$$
Observe that $\Delta_{\{v\},v}=D^2$.
For vertices $u,v$ of $K$ define {\bf double Whitehead product}
$[u,v]\in \pi_3(\mathcal Z_K)$ to be the homotopy class of the composition
$$
S^3 = \partial(D^2\times D^2) = \partial(\Delta_{\{u\},u}\times \Delta_{\{v\},v}) =
\partial\Delta_{\{u\},u}\times \Delta_{\{v\},v}\cup \Delta_{\{u\},u}\times \partial\Delta_{\{v\},v}
\subset \mathcal Z_K
$$
of the standard identifications and inclusion.\footnote{Observe that $[u,v]$ is not the Whitehead product of two elements in $\pi_2(\mathcal Z_K)=0$.}


We have $[u,v]=0$ if and only if $\{u,v\}$ is an edge of $K$.
Thus if $[u,v]=0$, then the composition
$$
D^4 = D^2\times D^2 = \Delta_{\{u,v\},u}\times \Delta_{\{u,v\},v} \subset \mathcal Z_K
$$
of the standard identifications and inclusion is a canonical inclusion extending to $D^4$ the above composition defining $u,v$.
Hence for vertices $u,v,w$ of $K$ such that $[u,v]=[v,w]=[w,u]=0$ we can define {\bf triple Whitehead product} $[u,v,w]\in \pi_5(\mathcal Z_K)$ to be the homotopy class of the composition
$$
S^5 = \partial(D^2\times D^2\times D^2) =
\partial(\Delta_{\{u\},u}\times \Delta_{\{v\},v} \times \Delta_{\{w\},w})\subset \mathcal Z_K
$$
of the standard identifications and the union of the above canonical
inclusions.

If all the $(r-1)$-fold Whitehead products of vertices $v_1,\ldots,v_r$ of $K$ are zero, then we can define
{\bf $r$-tuple Whitehead product} $[v_1,\ldots,v_r]\in \pi_{2r-1}(\mathcal Z_K)$ analogously, using a version of \cite[Proposition 3.3]{AP19} for $\mathcal Z_K$ instead of $(\C P^\infty)^K$ (this version is presumably proved analogously to \cite{AP19}).

An {\bf iterated higher Whitehead product} $w$ is a labelled tree.
To the expression
$$\left[v_1,v_2,\left[v_3,v_4,[v_5,v_6,v_7],v_8\right],\left[v_9,v_{10}\right]\right]$$
there corresponds the tree having a root, 4 vertices $v_1,v_2,\left[v_3,v_4,[v_5,v_6,v_7],v_8\right],[v_9,v_{10}]$ joined to the root, 6 vertices $v_3,v_4,[v_5,v_6,v_7],v_8,v_9,v_{10}$ on the second level and 3 vertices $v_5,v_6,v_7$ on the third level.
It is clear what are the edges.

In \cite[Construction 4.4]{AP19} one essentially introduces the complex $\partial\Delta_w$ associated with a labelled tree $w$.

According to \cite[two paragraphs before Definition 2.2]{AP19} one can define analogously to above
 {\bf evaluation of a labelled tree $w$ on a complex $K$}, denoted $w(K)\in\pi_n(\mathcal Z_K)$ for some $n$.
No definition of this evaluation (or of analogous evaluation in $\pi_n((\C P^\infty)^K)$) is given in \cite{AP19}.
Perhaps such a definition could be given if one states and proves certain analogue of \cite[Proposition 3.3]{AP19} for iterated products and for $\mathcal Z_K$ instead of $(\C P^\infty)^K$.
Assuming that such a definition could be given, \cite[Theorem 5.1]{AP19} is as follows.

\begin{theorem}\label{t:abrpan}
For any labelled tree $w$ we have $w(\partial\Delta_w)\ne0$.
\end{theorem}

\begin{remark}[Specific remarks on \cite{AP19}]\label{r:spec}
The statements of main results in \cite{AP19} use the notions of iterated higher Whitehead product
in $\pi_*(\mathcal Z_K)$, of realizability and of non-triviality.

The definition of iterated higher Whitehead product in $\pi_*(\mathcal Z_K)$ is not
given in \cite{AP19}.
In the 4th paragraph of p. 4 in \cite{AP19} (here and below I refer to numeration of the arxiv version 3)
it is written `... we consider {\it general iterated} higher Whitehead products, i.e. higher Whitehead products in which arguments can be higher Whitehead products'.
No definition of these objects considered is given.
It is not clear to me how to give the definition without stating and proving certain analogue of \cite[Proposition 3.3]{AP19} for iterated products.

The object used in Theorem 5.1 and in Definition 2.2 of realizability from \cite{AP19} is
{\it iterated higher Whitehead product in $\pi_*(\mathcal Z_K)$} not
{\it general iterated higher Whitehead product in $\pi_*((\C P^\infty)^K)$} attempted to be defined in the 4th paragraph of p. 4 in \cite{AP19}.
The phrase before Definition 2.2 does not define $w\in \pi_*(\mathcal Z_K)$
because the explanation `for dimensional reasons' of the existence of the lifting is unclear, because
it is not stated that the lifting is unique, and because no canonical way of choosing the lifting is described.
(This is resolved by the above direct working with $\mathcal Z_K$ and not mentioning the object $(\C P^\infty)^K$ which is not used in the statements of the main results.)



The definition of non-triviality used in \cite[Theorem 5.1]{AP19} is not given in \cite{AP19}.
The word `nontrivial' should simply be deleted from \cite[Theorem 5.1]{AP19}, cf. Theorem \ref{t:abrpan}. 


The first sentences of \cite[Theorem 5.1]{AP19} and of \cite[Construction 4.4]{AP19}) do not make sense because iterated higher Whitehead products are not defined in \cite{AP19} if $K$ is not specified.
(This is resolved by the above definitions of an iterated higher Whitehead product as a labelled tree, and of its evaluation on a complex.)
\end{remark}



\begin{thebibliography}{RSS}


\bibitem[AP19]{AP19} S. Abramyan and T. Panov, Higher Whitehead products in moment-angle complexes and substitutions of simplicial complexes, Tr. Mat. Inst. Steklova 305 (2019), Algebraicheskaya Topologiya Kombinatorika i Matematicheskaya Fizika, 7--28, arXiv:1901.07918v3.

\end{thebibliography}
\end{document}